\documentclass[12pt]{amsart}
\usepackage[foot]{amsaddr}
\usepackage[T1]{fontenc}
\usepackage[utf8]{inputenc}
\usepackage{color,graphicx,amsmath,amsfonts,amssymb,amsxtra,setspace,xspace,subfig, tikz, float}
\usepackage[style=numeric, citestyle=numeric-comp, backend=bibtex, sorting=none]{biblatex}
\usepackage{filecontents}
\usepackage{fancyhdr}
\fancypagestyle{goteastyle}{%
	\fancyhf{}%
	\fancyhead[OC]{\textsc{Mircea Gotea}}%
	\fancyhead[EC]{\textsc{An Extension of Pythagoras Theorem}}%
	\fancyfoot[C]{\thepage}%
}
\usetikzlibrary{arrows,shapes,calc}

\usepackage[colorlinks=true]{hyperref}
\hypersetup{urlcolor=blue, citecolor=red, linkcolor=blue}

\newtheorem{theorem}{Theorem}[]
\newtheorem*{lemma}{Lemma}

\newtheorem*{observation}{Observation}

\usepackage{etoolbox}
\makeatletter
\patchcmd{\@setauthors}{\MakeUppercase}{\large\scshape}{}{}
\patchcmd{\@setfoot@addresses}{\scshape}{}{}{}
\renewcommand{\@setemails}{\mbox{}$\,\,\,$ E-mail:\space{\ttfamily \emails}}
\makeatother

\newcommand{\PTLBL}[3][]{
\filldraw[white] (#2) circle (1.7pt);
\draw[black] (#2) circle (1.7pt) node[anchor=#3]{$#1#2$};
}

\DeclareFieldFormat[article]{title}{\emph{#1.\isdot}}
\begin{document}
\title[An Extension of Pythagoras theorem]
{An Extension of Pythagoras theorem}

\author{Mircea Gotea$^*$}
\email{mircea@gotea.ro}
\address{$^*$ Retired teacher of Mathematics, Gurghiu, Mure\c{s} County, Romania.}

\date{}

\begin{abstract}
 This article proves a Pythagoras-type formula for the sides and diagonals of a polygon inscribed in a semicircle having one of the sides of the polygon as diameter.
\end{abstract}

\maketitle


\par \vskip \baselineskip
\noindent 
{\large T}his paper presents a Pythagoras type theorem regarding a cyclic polygon having one of the sides as a diameter
of the circumscribed circle.\\
Let us first observe that Pythagoras theorem \cite{rusu, wagner, hadamard, bottema, saikia} and its reciprocal can be stated differently.
\begin{theorem}
 If a triangle is inscribed in a semicircle with radius $R$, such that its longest side is the diameter of the semicircle, then the following relation exists between its sides:
 \begin{eqnarray}
 a^2 & = & b^2 + c^2, \label{Pythagoras_theorem}
 \end{eqnarray}
which can also be expressed as
\begin{eqnarray}
4R^2 & = & b^2 + c^2.
\end{eqnarray}
Its reciprocal can be states as follows: if the relation $a^2 = b^2 + c^2$ exists between the sides of a triangle,
then that triangle can be inscribed in a semicircle such that the longest side is also the diameter of the semicircle.
\end{theorem}
We can state and prove a similar theorem for quadrilaterals. Let $ABCD$ be a quadrilateral with vertices $B$ and $C$ located on the semicircle that has the longest side $AD$ as diameter.

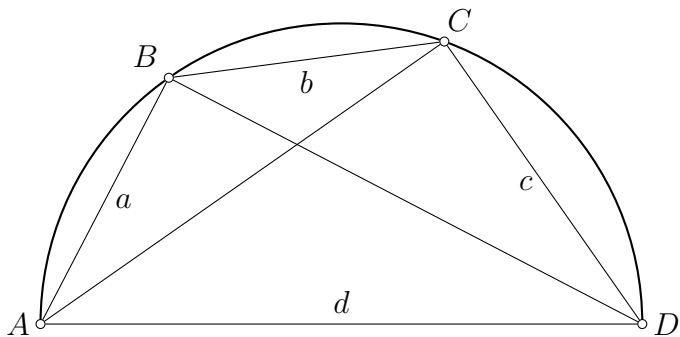
\begin{figure}[H]
\centering
{\center{\begin{tikzpicture}
\draw[thick] (4,0) arc (0:180:4);

\coordinate (A) at (-4,0);
\coordinate (D) at (4,0);
\coordinate (B) at ({4*cos(125)}, {4*sin(125)});
\coordinate (C) at ({4*cos(70)}, {4*sin(70)});

\draw[black] (A) -- (B);
\draw[black] (B) -- (C);
\draw[black] (C) -- (D);
\draw[black] (A) -- (C);
\draw[black] (B) -- (D);
\filldraw[black] ($ (A)!.5!(D) $) node[anchor=south]{$d$};
\filldraw[black] ($ (A)!.5!(B) $)  node[anchor=west]{$a$};
\filldraw[black] ($ (B)!.5!(C) $) node[anchor=north]{$b$};
\filldraw[black] ($ (C)!.5!(D) $) node[anchor=east]{$c$};

\draw[black] (A) -- (D);

\PTLBL{A}{east}
\PTLBL{B}{south east}
\PTLBL[\hspace{1em}]{C}{south}
\PTLBL{D}{west}

\end{tikzpicture}}}
\end{figure}

\begin{theorem}\label{quadrilateral_theorem}
If a quadrilateral $ABCD$ with its longest side $AD = d$ can be inscribed in a semicircle with the diameter $AD$, then its sides satisfy the following relation:
\begin{eqnarray}
d^2 & = & a^2 + b^2 + c^2 + \frac{2abc}{d}. \label{quadrilateral_Pythagoras}
\end{eqnarray}
\begin{proof}
We can successively deduce the following:
\begin{eqnarray*}
AD^2 & = & AB^2 + BD^2\\
& = & AB^2 + BC^2 + CD^2 - 2BC \cdot CD \cdot \cos(C)\\
& = & AB^2 + BC^2 + CD^2 + 2BC \cdot CD \cdot \cos(A)\\
& = & AB^2 + BC^2 + CD^2 + 2BC \cdot CD \cdot \frac{AB}{AD}.
\end{eqnarray*}
\end{proof}
\end{theorem}
\noindent {\bf Particular case.} If $B$ coincides with $A$, then $a = 0$, the quadrilateral becomes triangle $ACD$,
and relation (\ref{quadrilateral_Pythagoras}) becomes $d^2 = b^2 + c^2$, which represents
Pythagoras theorem for this right triangle.\\
The reciprocal of Theorem \ref{quadrilateral_theorem} does not hold true: if the sides of a quadrilateral
$ABCD$ satisfy relation (\ref{quadrilateral_Pythagoras}), the quadrilateral does not have to be inscribed in a semicircle having the longest side as diameter. For example, a quadrilateral with sides $AD = 4\sqrt{2},\ {AB = \sqrt{2}}$, ${BC = 3+\sqrt{5}},\, CD = 3 - \sqrt{5}$ and $m(\sphericalangle A) = 90^{\circ}$ does satisfy relation
(\ref{quadrilateral_Pythagoras}), has $AD$ as its longest side, but cannot be inscribed in a semicircle with diameter $AD$.\\
However, the following statement is true:
\begin{lemma}
Let four segments have the lengths $a$, $b$, $c$, and $d$. If they satisfy relation (\ref{quadrilateral_Pythagoras}),
then there is at least one quadrilateral with these sides that can be inscribed in a semicircle having the longest side as diameter.
\end{lemma}
\begin{proof}
On a circle with diameter $AD = d = 2R$, we set a point $B$ such that $AB = a$, and a point $C$ such that $BC = b$.
We can prove that $CD = c$ when the sides satisfy relation (\ref{quadrilateral_Pythagoras}).\\
We have successively:
\begin{eqnarray*}
CD & = & AD\cos(D) \\
& = & -AD\cos(B) \\
& = & -AD \cdot \frac{AB^2 + BC^2 - AC^2}{2 \cdot AB \cdot BC} \\
& = & -AD \cdot \frac{AB^2 + BC^2 - (AD^2 - CD^2)}{2 \cdot AB \cdot BC} \\
& = & -d \cdot \frac{a^2 + b^2 + c^2 - d^2}{2ab}   \\
& = & c .
\end{eqnarray*}

We thus found a quadrilateral with sides of lengths $a, b, c, d$.

\noindent Let us also observe that if $a, b$, and $c$ are distinct, there are three incongruent quadrilaterals that satisfy this requirement, depending on which of the three sides, $a, b$, or $c$, is selected as the opposite side to $d$.
\end{proof}
\begin{theorem}
If a pentagon $ABCDE$ with sides $AB = a$, $BC = b$, $CD = c$, $DE = d$, $AE = 2R$ is inscribed in a circle
with radius $R$, then its sides satisfy the following relation:
\begin{eqnarray}
4R^2 & = & a^2 + b^2 + c^2 + \frac{aby + xcd}{R}, \label{pentagon_Pythagoras}
\end{eqnarray}
where $x = AC$ and $y = CE$.
\end{theorem}
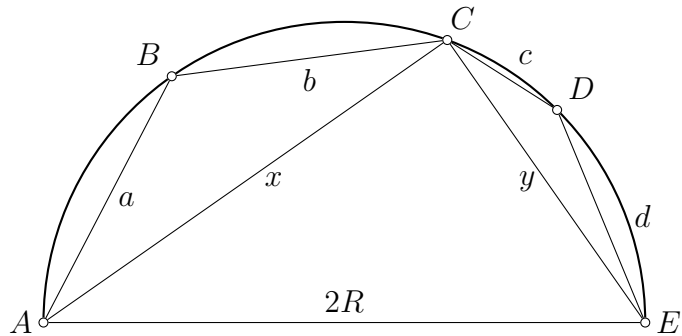
\begin{figure}[H]
\centering
{\center{\begin{tikzpicture}

\draw[thick] (4,0) arc (0:180:4);

\coordinate (A) at (-4,0);
\coordinate (B) at ({4*cos(125)}, {4*sin(125)});
\coordinate (C) at ({4*cos(70)}, {4*sin(70)});
\coordinate (D) at ({4*cos(45)}, {4*sin(45)});
\coordinate (E) at (4,0);

\draw[black] (A) -- (E);
\draw[black] (A) -- (B);
\draw[black] (B) -- (C);
\draw[black] (C) -- (D);
\draw[black] (D) -- (E);
\draw[black] (A) -- (C);
\draw[black] (C) -- (E);

\filldraw[black] ($ (A)!.5!(E) $) node[anchor=south]{$2R$};
\filldraw[black] ($ (A)!.5!(B) $)  node[anchor=west]{$a$};
\filldraw[black] ($ (B)!.5!(C) $) node[anchor=north]{$b$};
\filldraw[black] ($ (C)!.5!(D) $) node[anchor=south west]{$\,c$};
\filldraw[black] ($ (D)!.5!(E) + (0.3,0) $) node[anchor=west]{$d$};
\filldraw[black] ($ (C)!.5!(E) $) node[anchor=east]{$y$};
\filldraw[black] ($ (A)!.57!(C) $) node[anchor=north]{$x$};

\PTLBL{A}{east}
\PTLBL{B}{south east}
\PTLBL[\hspace{1em}]{C}{south}
\PTLBL{D}{south west}
\PTLBL{E}{west}
 
\end{tikzpicture}}}
\end{figure}

\begin{proof}
In the quadrilateral $ABCE$ we can apply relation (\ref{quadrilateral_Pythagoras}):
\begin{eqnarray}
4R^2 & = & a^2 + b^2 + y^2 + \frac{aby}{R} \label{temp}.
\end{eqnarray}
Using the Law of Cosines we have
\begin{eqnarray*}
y^2 & = & c^2 + d^2 - 2cd\cos(D) \nonumber\\
& = & c^2 + d^2 + 2cd\cos\left(\sphericalangle CAE\right) \nonumber\\
& = & c^2 + d^2 + 2cd \frac{x}{2R}.
\end{eqnarray*}
Replacing $y^2$ from this relation in (\ref{temp}), we obtain relation (\ref{pentagon_Pythagoras}).
\end{proof}

\noindent For a hexagon, we can similarly prove the following:

\begin{theorem}
If we inscribe a hexagon in a semicircle such that the longest side is also the diameter, then the following relation is true:
\begin{eqnarray}
4R^2 & = & a^2 + b^2 + c^2 + d^2 + e^2 + \frac{abz + ycx + ude}{R}.
\end{eqnarray}
\end{theorem}
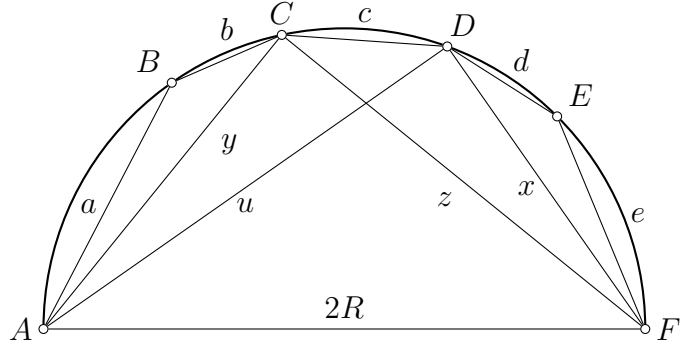
\begin{figure}[H]
\centering
{\center{\begin{tikzpicture}
\draw[thick] (4,0) arc (0:180:4);

\coordinate (A) at (-4,0);
\coordinate (B) at ({4*cos(125)}, {4*sin(125)});
\coordinate (C) at ({4*cos(102)}, {4*sin(102)});
\coordinate (D) at ({4*cos(70)}, {4*sin(70)});
\coordinate (E) at ({4*cos(45)}, {4*sin(45)});
\coordinate (F) at (4,0);

\draw[black] (A) -- (F);

\draw[black] (A) -- (B);
\draw[black] (B) -- (C);

\draw[black] (C) -- (D);
\draw[black] (D) -- (E);
\draw[black] (E) -- (F);

\draw[black] (A) -- (D);
\draw[black] (D) -- (F);

\draw[black] (A) -- (C);
\draw[black] (C) -- (F);

\filldraw[black] ($ (A)!.5!(F) $) node[anchor=south]{$2R$};
\filldraw[black] ($ (A)!.5!(B) $)  node[anchor=east]{$a$};
\filldraw[black] ($ (B)!.5!(C) + (0,0.1) $)  node[anchor=south]{$b$};
\filldraw[black] ($ (C)!.5!(D) + (0,0.1) $)  node[anchor=south]{$c$};

\filldraw[black] ($ (D)!.5!(E) $) node[anchor=south west]{$d$};
\filldraw[black] ($ (E)!.5!(F) + (0.25,0.1) $) node[anchor=west]{$e$};
\filldraw[black] ($ (D)!.5!(F) $) node[anchor=east]{$x$};
\filldraw[black] ($ (A)!.5!(D) $) node[anchor=north]{$u$};
\filldraw[black] ($ (A)!.7!(C) $) node[anchor=north west]{$y$};
\filldraw[black] ($ (C)!.5!(F) $) node[anchor=north east]{$z$};

\PTLBL{A}{east}
\PTLBL{B}{south east}
\PTLBL{C}{south}
\PTLBL[\hspace{1em}]{D}{south}
\PTLBL{E}{south west}
\PTLBL{F}{west}

\end{tikzpicture}}}
\end{figure}

\noindent We can restate the theorem in the case of a polygon with $n$ sides.
\begin{theorem}\label{general_theorem}
If a polygon $A_1A_2A_3 \dots A_{n-1}A_n$ can be inscribed in a semicircle with diameter $A_1A_n$, then the following relation is true:
\begin{eqnarray}
\left(A_1A_n\right)^2 & = & \sum_{k = 1}^{n - 1}\left(A_kA_{k + 1}\right)^2 +
2\sum_{k = 1}^{n - 3}\frac{\left(A_1A_{k + 1}\right)\left(A_{k + 1}A_{k + 2}\right)\left(A_{k + 2}A_n\right)}{A_1A_n}.
\label{general_formula}
\end{eqnarray}
\end{theorem}

\begin{observation}
The points $A_1$, $A_{k + 1}$, $A_{k + 2}$, and $A_n$, for $1 \leq k \leq n - 3$, involved in the sum
from the right side of formula (\ref{general_formula}), are the vertices of a cyclic quadrilateral
inscribed in a semicircle.
\end{observation}

\begin{proof}
By using the mathematical induction technique, we have already demonstrated the particular cases for
$n = 3$ and $n = 4$.\\
Let us assume that the property is true for any polygon with $n$ sides and let us consider the case of a polygon $A_1A_2A_3 \dots A_nA_{n+1}$, for $n \geq 4$, inscribed in a circle with diameter $A_1A_{n+1}$.
Following the induction hypothesis, in an $n$-sided polygon $A_1A_2A_3 \dots A_{n - 1}A_{n+1}$ the following relation is true:
	\begin{eqnarray}\label{proof_of_general}
\left(A_1A_{n + 1}\right)^2 & = & \sum_{k = 1}^{n - 2}\left(A_kA_{k + 1}\right)^2 + \left(A_{n - 1}A_{n + 1}\right)^2\\
	& & + 2\sum_{k = 1}^{n - 3}\frac{\left(A_1A_{k + 1}\right)\left(A_{k + 1}A_{k + 2}\right)\left(A_{k + 2}A_{n + 1}\right)}{A_1A_{n + 1}}. \nonumber
\end{eqnarray}
Applying the Law of Cosines in triangle $A_{n - 1}A_nA_{n + 1}$, and using the facts that:
\begin{eqnarray*}
m\left(\sphericalangle A_{n - 1}A_nA_{n + 1}\right) + m\left(\sphericalangle A_{n - 1}A_1A_{n + 1}\right) & = &
180^{\circ}
\end{eqnarray*}
and
\begin{eqnarray*}
m\left(\sphericalangle A_1A_{n - 1}A_{n + 1}\right) & = & 90^{\circ},
\end{eqnarray*}
we obtain:
\begin{eqnarray*}
\left(A_{n - 1}A_{n + 1}\right)^2 & = & \left(A_{n - 1}A_n\right)^2 + \left(A_nA_{n + 1}\right)^2
- 2\left(A_{n - 1}A_n\right)\left(A_nA_{n + 1}\right)\cos\left(A_n\right)\\
& = & \left(A_{n - 1}A_n\right)^2 + \left(A_nA_{n + 1}\right)^2
+ 2\left(A_{n - 1}A_n\right)\left(A_nA_{n + 1}\right)\cos\left(\sphericalangle A_{n - 1}A_1A_{n + 1}\right)\\
& = & \left(A_{n - 1}A_n\right)^2 + \left(A_nA_{n + 1}\right)^2
+ 2\left(A_{n - 1}A_n\right)\left(A_nA_{n + 1}\right)\cdot\frac{A_1A_{n - 1}}{A_1A_{n + 1}}.
\end{eqnarray*}
	Substituting $(A_{n - 1}A_{n + 1})^2$ into relation (\ref{proof_of_general}), we obtain:
\begin{eqnarray*}
\left(A_1A_{n + 1}\right)^2 & = & \sum_{k = 1}^{n - 2}\left(A_kA_{k + 1}\right)^2 +
\left(A_{n - 1}A_n\right)^2 + \left(A_nA_{n + 1}\right)^2 + 2\left(A_{n - 1}A_n\right)\left(A_nA_{n + 1}\right)\cdot\frac{A_1A_{n - 1}}{A_1A_{n + 1}}\\
& \ &
+ 2\sum_{k = 1}^{n - 3}\frac{\left(A_1A_{k + 1}\right)\left(A_{k + 1}A_{k + 2}\right)\left(A_{k + 2}A_{n + 1}\right)}{A_1A_{n + 1}}\\
& = & \sum_{k = 1}^n\left(A_kA_{k + 1}\right)^2
+ 2\sum_{k = 1}^{n - 2}\frac{\left(A_1A_{k + 1}\right)\left(A_{k + 1}A_{k + 2}\right)\left(A_{k + 2}A_{n + 1}\right)}{A_1A_{n + 1}}.
\end{eqnarray*}
The proof is now complete.
\end{proof}
\newpage
We can also note the validity of a reciprocal of this theorem:
	Let $A_1A_2$, $A_2A_3$, \dots, $A_{n-1}A_n$, $A_nA_1$ be $n$ segments with lengths $A_1A_2 = a_1$, $A_2A_3 = a_2$,
\dots, $A_{n-1}A_n = a_{n-1}$, $A_nA_1 = a_n$. If these segments satisfy relation (\ref{general_formula}), then there is at least one polygon with these segments as sides that can be inscribed in the semicircle with diameter $A_1A_n$.
	The proof can be developed in the same manner as for the case of $n = 4$.\\

{\bf Acknowledgments:} \ I would like to thank Professor Aurel I. Stan, Department of Mathematics, The Ohio State University, U.S.A., for his kind advice and encouragements.\\

{\bf Note:} \ This article was originally published in Romanian here \cite{gotea}.

\nocite{*}
\printbibliography
\end{document}